\documentclass{article}
\usepackage{amsmath, amsthm, amscd, amsfonts, amssymb, graphicx, color}

\newtheorem{theorem}{Theorem}[section]

\newtheorem{proposition}[theorem]{Proposition}
\newtheorem{corollary}[theorem]{Corollary}

\def\title#1{{\Large\bf  \begin{center} #1 \vspace{0pt} \end{center}  } }
\def\authors#1{{\large\bf \begin{center} #1 \vspace{0pt} \end{center} } }
\def\university#1{{\sl \begin{center} #1 \vspace{0pt} \end{center} } }
\def\inst#1{\unskip $^{#1}$}
\usepackage{amsfonts}

\begin{document}

%
%
%

\title{Inertia of Loewner Matrices}

\bigskip

%
%

\authors{Rajendra Bhatia\inst{1}, 
Shmuel Friedland\inst{2}, 
Tanvi Jain\inst{3} 
}


\smallskip

%
%

\university{\inst{1} Indian Statistical Institute, New Delhi 110016, India\\rbh@isid.ac.in}
\university{\inst{2} Department of Mathematics, Statistics and Computer Science, University of Illinois at Chicago, Chicago, 60607-7045, USA \\ friedlan@uic.edu}
\university{\inst{3} Indian Statistical Institute, New Delhi 110016, India\\tanvi@isid.ac.in}%

\begin{abstract}
Given positive numbers $p_1 < p_2 < \cdots < p_n,$ and a real number $r$ let 
$L_r$ be the $n \times n$ matrix with its $i,j$ entry equal to 
$(p_i^r-p_j^r)/(p_i-p_j).$ A 
well-known theorem of C. Loewner says that $L_r$ is positive definite when $0 < 
r < 1.$ In contrast, R. Bhatia and J. Holbrook, (Indiana Univ. Math. J, 49 
(2000) 1153-1173) showed that when $1 < r < 2,$ the matrix $L_r$ has only one 
positive eigenvalue, and made a conjecture about the signatures of eigenvalues 
of $L_r$ for other $r.$ That conjecture is proved in this paper.
\end{abstract}



\section{Introduction} 
Let $f$ be a real-valued $C^{1}$ function on $(0, \infty).$ Let $p_1 < p_2 < 
\cdots < p_n$ be any $n$ points in $(0, \infty).$ The $n \times n$ matrix
\begin{equation}
L_f (p_1, \ldots, p_n) = \left [\frac{f(p_i) - f (p_j)}{p_i - p_j} 
\right ]_{i,j=1}^n  \label{eq1}
\end{equation}
is called a {\it Loewner matrix} associated with $f.$ It is understood that 
when $i=j,$ the quotient in (\ref{eq1}) represents the limiting value 
$f^{\prime} (p_i).$ Of particular interest to us are the functions         
$f(t) = t^r,$ $r \in \mathbb{R}.$ In this case we write $L_r$ for $L_f 
(p_1, \ldots, p_n),$  where the roles of $n$ and $p_1, \ldots p_n$  can be 
inferred from the context. Thus $L_r$ is the $n \times n$ matrix
\begin{equation}
L_r = \left [ \frac{p_i^r - p_j^r}{p_i-p_j}\right ]_{i,j=1}^n.   \label{eq2}
\end{equation}

Loewner matrices are important in several contexts, of which we mention two 
that led to the present study. (The reader may see Section 4.1 of \cite{hn} for 
an excellent discussion of both these aspects of Loewner matrices.)  The 
function $f$ on $(0, \infty)$ induces, via 
the usual functional calculus, a matrix function $f(A)$ on the space of 
positive definite matrices. Let $Df(A)$ be the Fr\'echet derivative of this 
function.This is a linear map on the space of Hermitian matrices. The {\it 
Daleckii-Krein formula} describes the action of this map in terms of Loewner 
matrices. Choose an orthonormal basis in which $A=\mbox{diag}(p_1, 
\ldots, p_n).$ Then the formula says that for every Hermitian $X$

\begin{equation}
Df(A) (X) = L_f (p_1, \ldots, p_n) \circ X,   \label{eq3}
\end{equation}
where $A \circ B$ stands for the entrywise product $[a_{ij} b_{ij}]$ of $A$ and 
$B.$

The function $f$ is said to be {\it operator monotone} on $(0, \infty)$ if $A 
\ge B > 0$ implies $f(A) \ge f(B).$ (As usual $A \ge 0$ means $A$ is positive 
semidefinite.) A fundamental theorem due to Charles Loewner says that $f$ is 
operator monotone if and only if all Loewner matrices associated with $f$ (for 
every $n$ and for every choice $p_1, \ldots, p_n$) are positive semidefinite. 
Another basic fact, again proved first by Loewner, says that $f(t) = t^r$ is 
operator monotone if and only if $0 \leq r \leq 1.$ See \cite{rbh} Chapter V.

Combining these various facts with some well-kown theorems on positive linear 
maps \cite{rbh1} one can see that if $f$ is operator monotone, then the norm of 
$Df(A)$ obeys the relations
\begin{equation}
\|Df(A)\| = \|Df (A) (I)\| = \| f^{\prime} (A) \|,    \label{eq4}
\end{equation}
and is therefore readily computable. In particular, for the function $f(t) = 
t^r$ if we write $DA^r$ for $Df(A),$ then (\ref{eq4}) gives
\begin{equation}
\|DA^r\| = \| r A^{r-1}\|, \quad \mbox{for}\,\, 0 \leq r \leq 1.   \label{eq5}
\end{equation}
This was first noted in \cite{rb}, and used to derive perturbation bounds for the 
operator absolute value. Then in \cite{bs} Bhatia and Sinha showed that the 
relation (\ref{eq5}) holds also for $- \infty < r < 0$ and for $2 \leq r < 
\infty$ but, mysteriously, not for $1 < r < \sqrt{2}.$ The case $\sqrt{2} \leq 
r < 2,$ left open in this paper, was resolved in \cite{bhol} by Bhatia and 
Holbrook, who showed that here again the relation (\ref{eq5}) is valid.

One ingredient of the proof in \cite{bhol} is their Proposition 2.1 which says 
that when $1 < r < 2,$ the $n \times n$ matrix $L_r$ has just one positive 
eigenvalue. We have remarked earlier that when $0 < r < 1,$ the matrix $L_r$ is 
positive semidefinite and therefore, none of its eigenvalues is negative. This 
contrast as $r$ moves from $(0,1)$ to $(1,2)$ is intriguing, and raises the 
natural question about the behaviour of eigenvalues of $L_r$ for other values 
of $r.$ Bhatia and Holbrook \cite{bhol} made a conjecture about this behaviour and 
established a small part of it: they settled the cases $r=1,2,\ldots, n-1$ 
apart from $0 < r 
< 1$ and $1 < r < 2$ already mentioned. The main goal of this paper is to prove 
this conjecture in full. This is our Theorem 1.1.

Let $A$ be an $n \times n$ Hermitian matrix. The {\it inertia} of $A$ is the 
triple
$${\rm In} (A) = (\pi(A), \zeta (A), \nu(A)),$$
where $\pi(A)$ is the number of positive eigenvalues of $A,$ $\zeta (A)$ is the 
number of zero eigenvalues of $A,$ and $\nu(A)$ the number of negative 
eigenvalues of $A.$ Theorem \ref{thm1.1} describes the inertia of $L_r$ as $r$ 
varies over $\mathbb{R}.$ It is easy to see that the inertia of $L_{-r}$ is the 
opposite of the inertia of $L_r;$ i.e. $\pi (L_{-r}) = \nu (L_r)$ and $\nu 
(L_{-r}) = \pi (L_r).$ So we confine ourselves to the case $r > 0.$

\begin{theorem}\label{thm1.1}
Let $p_1 < p_2 < \cdots < p_n$ and $r$  be any positive real numbers and let 
$L_r$ be the matrix defined in (\ref{eq2}). Then
\begin{itemize}
\item[{\rm (i)}] $L_r$ is singular if and only if $r=1, 2, \ldots, n-1.$

\item[{\rm (ii)}] At the points $r=1,2,\ldots, n,$ the inertia of $L_r$ is 
given as follows:
$$r=2k \Rightarrow {\rm In}\,\, (L_r)\,\, = \,\,(k, n-r, k), $$
and
$$r=2k-1  \Rightarrow {\rm In}\,\, (L_r) \,\,=\,\, (k, n-r, k-1).$$

\item[{\rm (iii)}] If $0 < r < n$ and $r$ is not an integer,  then
$$\lfloor r \rfloor =2k \Rightarrow {\rm In}\,\, (L_r) \,\,= \,\,(n-k, 0, k) $$
and
$$\lfloor r \rfloor =2k-1  \Rightarrow {\rm In}\,\, (L_r) \,\,= \,\,(k, 0, 
n-k).$$

\item[{\rm (iv)}] If $r > n-1,$  \,\, then ${\rm In}\,\, (L_r) = {\rm In} 
\,\, (L_n).$

\item[{\rm (v)}] Every nonzero eigenvalue of $L_r$ is simple.
\end{itemize}
\end{theorem}

It is helpful to illustrate the theorem by a picture. Figure 1 is a diagram of the (scaled) eigenvalues of a $6 \times 6$ matrix $L_r$ when $p_i$ are fixed 
and $r$ varies. Some of the eigenvalues are very close to zero. To be able to distinguish between them the vertical scale has been expanded. 

\vskip0.1in
\begin{center}
\begin{figure}[ht]
\includegraphics[width=4.5in]{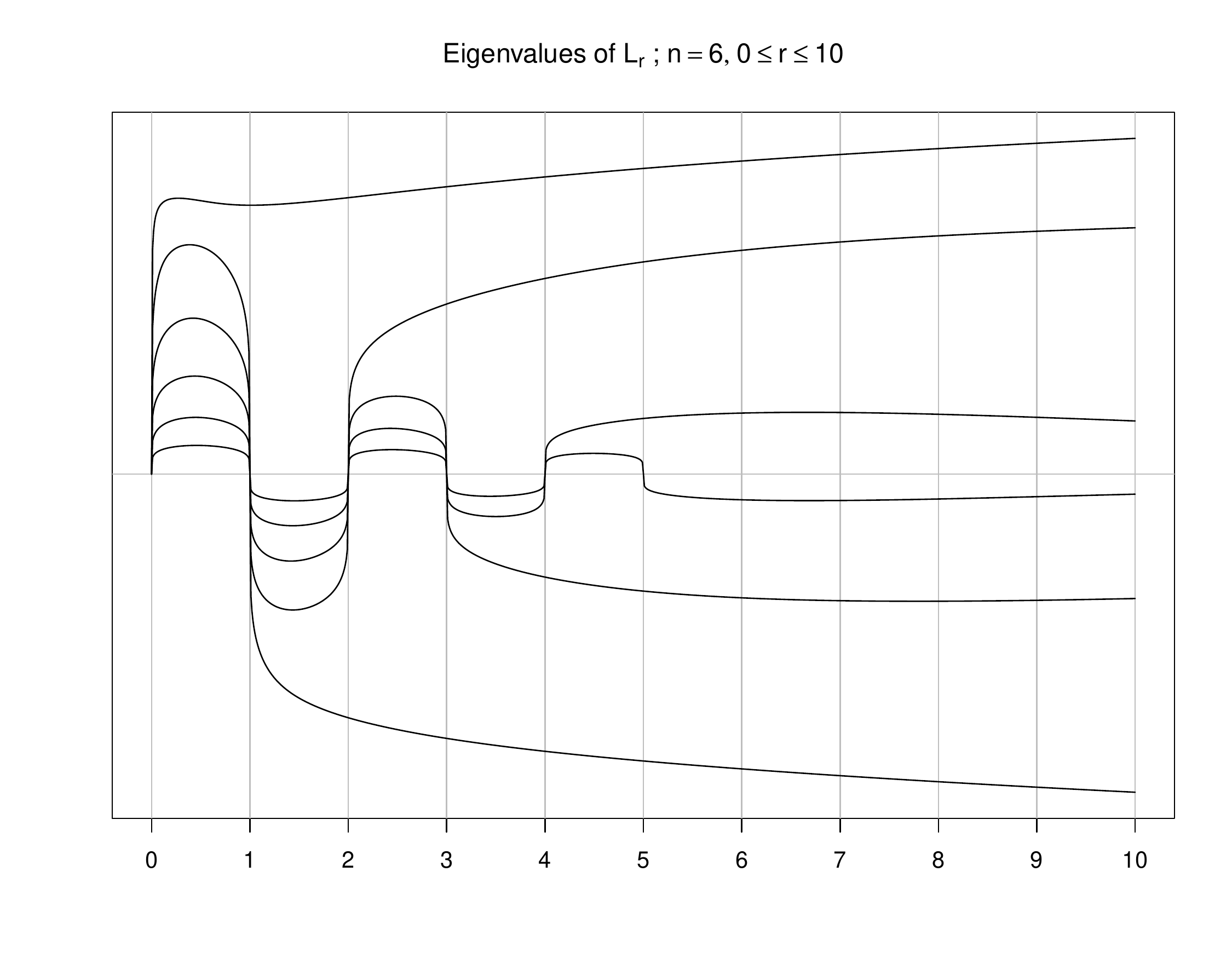}
\caption{}
 \end{figure}
\end{center}
\vskip0.1in

We have already mentioned that for $0 < r < 1,$ statement (iii) of Theorem 
\ref{thm1.1} follows from Loewner's theorem, and for $1 < r < 2$ it was 
established in \cite{bhol}. The case $2 < r < 3$ was accomplished by Bhatia and 
Sano in \cite{rbsano}. We briefly explain this work.

Let $\mathcal{H}_1$ be the space
\begin{equation}
 \mathcal{H}_1 = \left \{ x = (x_1, \ldots, x_n) : \sum_{i=1}^{n} x_i = 0 
\right \}.  \label{eq6}
\end{equation}
An $n \times n$ Hermitian matrix $A$ is said to be {\it conditionally positive 
definite} if $\langle x, Ax \rangle \ge 0$ for all $x \in\mathcal{H}_1,$ and if 
$-A$ has this property, then we say that $A$ is {\it conditionally negative 
definite}. Since $\mbox{dim} \mathcal{H}_1 = n-1,$ a nonsingular conditionally 
positive definite matrix which is not positive definite has inertia $(n-1, 0, 
1).$

In \cite{rbsano} it was shown that when $1 < r < 2,$ the matrix $L_r$ is nonsingular 
and conditionally negative definite. It follows that ${\rm In} \,\, (L_r) = 
(1,0, n-1),$ 
a fact established earlier in \cite{bhol}. It was also shown in \cite{rbsano} that when 
$2 < r < 3,$ the matrix $L_r$ is nonsingular and conditionally positive 
definite. From this it follows that ${\rm In} \,\,(L_r) = (n-1, 0 , 1).$

More generally, Bhatia and Sano \cite{rbsano} showed that $f$ on $(0, \infty)$ is {\it 
operator convex} if and only if all Loewner matrices $L_f$ are conditionally 
negative definite. This is a characterisation analogous to Loewner's for 
operator monotone functions. It is well-known that $f(t) = t^r$ is operator 
convex for $1 \leq r \leq 2.$

The proof of Theorem \ref{thm1.1} is given in Section 2. We also indicate how 
the proofs for the parts already given in \cite{bhol} and \cite{rbsano} can be considerably 
simplified. The inertia of the matrix $[(p_i+p_j)^r]$ has been studied by 
Bhatia and Jain in \cite{bj}. Some ideas in our proofs are similar to the ones 
used there.

\section{Proofs and Remarks}

Let $X$ be an $n \times n$ nonsingular matrix. The transformation $A \mapsto 
X^{\ast} AX$ on Hermitian matrices is called a {\it congruence.} The {\it 
Sylvester Law of Inertia} says that
\begin{equation}
{\rm In} \,\, (X^{\ast} AX) = {\rm In} \,\, A \,\,\,\mbox{for all}\,\,\, X 
\in GL(n).   \label{eq7}
\end{equation}
Let $D$ be the diagonal matrix
\begin{equation}
D = \mbox{diag}\, (p_1, \ldots, p_n).   \label{eq8}
\end{equation}
Then for every $r$
\begin{equation}
L_{-r} = - D^{-r} L_r D^{-r}.   \label{eq9}
\end{equation}
Hence by Sylvester's Law
\begin{equation}
{\rm In} \,\, L_r = (i_1, i_2, i_3) \Leftrightarrow {\rm In} \,\, L_{-r} = 
(i_3, i_2, i_1).   
\label{eq10}
\end{equation}
Thus all statements about ${\rm In} \,\, L_r$ for $r > 0$ give information 
about ${\rm In} \,\,L_{-r}$ as well.

Make the substitution $p_i = e^{2x_{i}},$ $x_i \in \mathbb{R}.$ A simple 
calculation shows that
$$L_r = \left [\frac{e^{rx_{i}}}{e^{x_{i}}} \,\, \frac{\sinh \,r 
(x_{i}-x_{j})}{\sinh (x_{i}-x_{j})}  \,\, \frac{e^{rx_{j}}}{e^{x_{j}}}\right ]. 
$$
In other words,
\begin{equation}
L_r = \Delta \widetilde{L}_r  \Delta,  \label{eq11}
\end{equation}
where $\Delta = \mbox{diag}\,(e^{(r-1)x_{1}}, \ldots, e^{(r-1)x_{n}} ),$ and 
\begin{equation}
 \widetilde{L}_r = \left [\frac{\sinh \,r 
(x_{i}-x_{j})}{\sinh (x_{i}-x_{j})}   \right ].\label{eq12}
\end{equation}
By Sylvester's Law ${\rm In} \,\, L_r = {\rm In} \,\, \widetilde{L}_r.$ 
Several properties of $L_r$ 
can be studied via $\widetilde{L}_r,$ and vice versa. This has been a very 
effective tool in deriving operator inequalities; see, the work of Bhatia and 
Parthasarathy \cite{bp} and that of Hiai and Kosaki \cite{hk1, hk2, hk3, hk4}.

When $n=2$ we have
$$\widetilde{L}_r = \left [\begin{array}{cc} r &  \frac{\sinh \,r 
(x_{1}-x_{2})}{\sinh (x_{1}-x_{2})}  \\  \frac{\sinh \,r 
(x_{1}-x_{2})}{\sinh (x_{1}-x_{2})}  & r \end{array}  \right ]. $$
So $\det \widetilde{L}_r = r^2 - \frac{\sinh^2 \,r 
(x_{1}-x_{2})}{\sinh^2 (x_{1}-x_{2})}.$ Thus  $\det \widetilde{L}_r$ is 
positive for $0 < r < 1,$ zero for $r=1,$ and negative for $r>1.$ One 
eigenvalue of $\widetilde{L}_r$ is always positive, and this shows that the 
second eigenvalue is positive, zero, or negative depending on whether $0< r < 
1,$ $r=1,$ or $r>1,$ respectively. This establishes Theorem \ref{thm1.1} in the 
simplest case $n=2.$  

An interesting corollary can be deduced at this stage. According to the two 
theorems  of Loewner mentioned in Section 1, $f$ is operator monotone if and 
only if all Loewner matrices $L_f$ are positive semidefinite, and $f(t) = t^r$ 
is operator monotone if and only if $0 \leq r \leq 1.$ Consequently, if $r>1,$ 
then there exists an $n,$ and positive numbers $p_1,\ldots,p_n$ such that the 
associated Loewner matrix (\ref{eq2}) is not positive definite. We can assert 
more:

\begin{proposition}\label{prop2.1}
Let $r > 1.$ Then for every $n \ge 2,$ and for every choice of $p_1, \ldots, 
p_n,$ the matrix $L_r$ defined in (\ref{eq2}) has at least one negative 
eigenvalue.
\end{proposition}

\vskip0.1in
\noindent{\it Proof}\,\,\, Consider  the $2 \times 2$ top left submatrix of 
$L_r.$ This is a Loewner matrix. By Theorem \ref{thm1.1} it has one negative 
eigenvalue. So, by Cauchy's interlacing principle, the $n \times n$ matrix 
$L_r$ has at least one negative eigenvalue. $\hfill{\square}$
\vskip0.1in
The Sylvester Law has a generalisation that is useful for us. Let $n \ge r,$ 
and let 
$A$ be an $r \times r$ Hermitian matrix and $X$ an $r \times n$ matrix of rank 
$r.$ Then

\begin{equation}
{\rm In} \,\, X^{\ast} AX = {\rm In} \,\, A + (0, n-r, 0).   \label{eq13}
\end{equation}
A proof of this may be found in \cite{bj}. This permits a simple transparent proof 
of Part 
(ii) of Theorem \ref{thm1.1}. (This part has already been proved in \cite{bhol}.) 
When $r$ is a positive integer we have
$$L_r = \left [p_i^{r-1} + p_i^{r-2} p_j + \cdots + p_j^{r-1} \right ] = 
W^{\ast} VW, $$
where $W$ is the $r \times n$ Vandermonde matrix 
$$W = \left [\begin{array}{cccc} 1 & 1 & \ldots & 1 \\ p_1 & p_2 & \ldots & p_n 
\\ \cdot & \cdot & \cdots & \cdot \\ p_1^{r-1} & p_2^{r-1} & \cdots & 
p_n^{r-1} \end{array} \right ], $$
and $V$ is the $r \times r$ antidiagonal matrix with all entries $1$ on its 
sinister diagonal and all its other entries equal to $0.$ If $r = 2 k,$ the 
matrix $V$ has $k$ of 
its eigenvalues equal to $1,$ and the other $k$ equal to $-1.$ If $r = 2k-1,$ 
then $k$ of its eigenvalues are equal to $1,$ and $k-1$ are equal to $-1.$ So, 
statement (ii) of Theorem \ref{thm1.1} follows from the generalised Sylvester's 
Law (\ref{eq13}). Next we prove statement (i).

Let $c_1, c_2, \ldots, c_n$ be real numbers, not all of which are zero. Let $f$ 
be the function on $(0, \infty)$ defined as
\begin{equation}
f(x) = \sum_{j=1}^{n} c_j \frac{x^r - p_j^r}{x - p_j}.   \label{eq14}
\end{equation}

\begin{theorem}\label{thm2.2}
Let $r$ be a positive real number not equal to $1,2,\ldots, n-1.$ Then the 
function 
$f$ defined in (\ref{eq14}) has at most $n-1$ zeros in $(0, \infty).$
\end{theorem}

\vskip0.1in
\noindent{\it Proof}\,\,\, Let $r_1 < r_2 < \cdots < r_m,$ and let $a_1, 
\ldots, a_m$ be real numbers not all of which are zero. Then the function 
\begin{equation}
g(x) = \sum_{j=1}^{m} a_j x^{r_{j}},   \label{eq15}
\end{equation}
has at most $m-1$ zeros in $(0,\infty).$ This is a well-known fact, and can be found in e.g., \cite{gg}, p.46.

Now let $f$ be the function defined in (\ref{eq14}) and let
\begin{equation}
g(x) = f(x) \prod_{j=1}^{n} (x-p_j).   \label{eq16}
\end{equation}
Then $g$ can be expressed in the form (\ref{eq15}) with $m=2n$ and
$$\left \{ r_1, \ldots, r_{2n} \right \} = \left \{ 0,1, \ldots, n-1, r, r+1, 
\ldots, r+n-1\right \}. $$
Further, we have $g(x)=x^rh_1(x)-h_2(x),$ where
\begin{align*}
h_1(x)=\sum\limits_{i=1}^{n}c_i\prod\limits_{j\ne i}(x-p_j),\ h_2(x)=\sum\limits_{i=1}^{n}c_ip_i^r\prod\limits_{j\ne i}(x-p_j).
\end{align*}
Both $h_1$ and $h_2$ are Lagrange interpolation polynomials of degree at most $n-1$. Since not all $c_i$ are zero, neither of these polynomials is identically zero. So, if $r\ne 1,2,\ldots,n-1,$ then $g$ is not the zero function.

Hence the function $g$ defined by (\ref{eq16}) has at most $2n-1$ zeros in 
$(0, \infty).$ Of these,  $n$ zeros occur at $x = p_j,$ $1 \leq j \leq n.$ So 
$f$ 
has at most $n-1$ zeros in $(0, \infty).$ $\hfill{\square}$

\begin{corollary}\label{cor2.3}
 Let $r$ be a  positive real number different from $1,2, \ldots, n-1.$ Then the 
matrix $L_r$ defined in (\ref{eq2}) is nonsingular. 
\end{corollary}

\vskip0.1in
\noindent{\it Proof}\,\,\, The matrix $L_r$ is singular if and only if there 
exists a nonzero vector $c=(c_1, \ldots, c_n)$ such that $L_r (c) = 0.$ In 
other words there exist real numbers $c_1, \ldots, c_n,$ not all zero, such 
that
$$\sum_{j=1}^{n} c_j \frac{p_i^r - p_j^r}{p_i - p_j} = 0 $$
for $i = 1,2, \ldots, n.$ But then the function $f(x)$ in (\ref{eq14}) would 
have $n$ zeros, viz., $x=p_1, \ldots, p_n.$ That is not 
possible.$\hfill{\square}$

We have proved Part (i) of Theorem \ref{thm1.1}. Part (iv) follows from this. 
If the inertia of $L_r$ were to change at some point $r_0 > n-1,$ then one of 
the eigenvalues has to change sign at $r_0.$ This is ruled out as $L_r$ is 
nonsingular for all $r>n-1.$

\vskip0.1in
Our argument shows that if $p_1 < p_2 < \cdots < p_n$ and $q_1 < q_2 < \cdots < 
q_n$ are two $n$-tuples of positive real numbers, then the matrix $\left 
[\frac{p_i^r - q_j^r}{p_i - q_j} \right ]$ is nonsingular for every positive 
$r$ different from $1,2, \ldots, n-1.$

An  $n \times n$ real matrix $A$ is said to be {\it strictly sign-regular} (SSR 
for short) if for every $1 \leq k \leq n,$ all $k \times k$ sub-determinants of 
$A$ are nonzero and have the same sign. If this is true for every $1 \leq k 
\leq r$ for some $r<n,$ then we say that $A$ is in the class SSR$_{\mbox{r}}.$ 
Sign-regular matrices and kernels are studied extensively in \cite{k}.

We have noted above that if $r$ is any positive real number and $k$ is any 
positive integer not greater than $r,$ then every $k \times k$ matrix of the 
form $\left [\frac{p_i^r - q_j^r}{p_i - q_j} \right ]$ is nonsingular. Let 
$L_r$ be an $n \times n$ Loewner matrix. Let $r \neq 1,2,\ldots, n-1.$ Using a 
homotopy argument one can see that all $k \times k$ sub-determinants of $L_r$ 
are nonzero and have the same sign. Thus $L_r$ is an SSR matrix. If 
$r=1,2,\ldots, n-1,$ then the same argument shows that for $k \leq r$  all 
$k \times k$ sub-determinants of $L_r$ are nonzero and have the same sign. In 
other words, $L_r$ is an  SSR$_{\mbox{r}}$ matrix.

Let $A$ be any matrix with eigenvalues $\lambda_1,\lambda_2,\cdots,\lambda_n$ arranged so that $|\lambda_1|\ge|\lambda_2|\ge \cdots\ge |\lambda_n|$. The Perron theorem tells us that if $A$ is entrywise positive, then $\lambda_1>0$ and $\lambda_1$ is a simple eigenvalue of $A$. (See\cite{hj}, p. 526). Applying this to successive exterior 
powers, we see that all eigenvalues of an SSR matrix are simple, and the $r$ 
nonzero eigenvalues of an SSR$_{\mbox{r}}$ matrix of rank $r$ are simple. This 
proves Part (v) of Theorem \ref{thm1.1}.

We now turn to proving Part (iii). Using the identity
$$\frac{p_i^r - p_j^r}{p_i - p_j} = \frac{p_i^{r-1} (p_i-p_j) + p_i 
(p_i^{r-2}-p_j^{r-2}) p_j + (p_i - p_j)p_j^{r-1}  }{p_i-p_j} $$
we see that for every $r \in \mathbb{R},$
\begin{equation}
L_r = D^{r-1} E + D L_{r-2} D + E D^{r-1},   \label{eq18}
\end{equation}
where $D$ is the diagonal matrix in (\ref{eq8}) and $E$ is the $n \times n$ 
matrix with all its entries equal to one.

By Loewner's Theorem $L_r$ is positive definite for $0 < r < 1,$ and because of 
(\ref{eq10}) it is negative definite for $-1 < r < 0.$ Now suppose $1 < r < 2.$ 
Let $x$ be any nonzero vector in the space $\mathcal{H}_1$ defined in 
(\ref{eq6}). Note that this $(n-1)$-dimensional space is the kernel of the 
matrix $E.$ Using (\ref{eq18}) we have
$$\langle x, L_r x \rangle = \langle x, D^{r-1} Ex \rangle + \langle x,  D 
L_{r-2} Dx \rangle + \langle x, ED^{r-1} x \rangle. $$
The first and the third term on the right hand side are zero because $Ex=0.$ So,
$$\langle x, L_r x \rangle  = \langle y, L_{r-2} y \rangle,$$
where $y=Dx.$ The last inner product is negative because $L_{r-2} < 0.$ Thus 
$\langle x, L_{r} x \rangle < 0$ for all $x \in \mathcal{H}_1.$  In other 
words, $L_r$ is conditionally negative definite if $1 < r < 2.$ The same 
argument shows that $L_r$ is conditionally positive definite if $2 < r < 3$ 
(because in this case $L_{r-2}$ is positive definite). This was proved in 
\cite{rbsano} by more elaborate arguments. In particular, we have

\begin{equation}
{\rm In} \,\, L_r = (1,0,n-1), \,\,\,\mbox{if}\,\,\, 1 < r < 2,   
\label{eq19}
\end{equation}
and
\begin{equation}
 {\rm In} \,\, L_r = (n-1, 0, 1), \,\,\,\mbox{if}\,\,\, 2 < r < 3.   
\label{eq20}
\end{equation}
We note here that if $n=3,$ then because of Part (iv) already proved we have 
${\rm In} \,\,L_r = (2,0,1)$ for all $r>2.$ So the theorem is completely 
proved for $n=3.$

Let $n > 3$ and suppose $3 < r < 4.$ Now consider the space
\begin{eqnarray*}
 \mathcal{H}_2 &=& \left \{ x : \sum x_i = 0, \,\, \sum p_i x_i =0 \right \}   
\\
 &=& \left \{ x : E x = 0, \,\,EDx=0 \right \}.
\end{eqnarray*}
This space is of dimension $n-2,$ being the orthogonal complement of the span 
of the vectors $e = (1,1, \ldots 1)$ and $p=(p_1, p_2, \ldots, p_n).$ Let $x 
\in \mathcal{H}_2.$ Again using the relation (\ref{eq18}) we see that
$$ \langle x, L_r x \rangle = \langle y, L_{r-2} y \rangle,$$
where $y=Dx.$ Since $EDx = 0,$ $y$ is in $\mathcal{H}_1,$ and since $1<r-2<2,$ 
we have $\langle x, L_r  x \rangle < 0.$ This is true for all $x \in 
\mathcal{H}_2.$ So, by the minmax principle $L_r$ has at least $n-2$ negative 
eigenvalues. The case $n=3$ of the theorem already proved shows that $L_r$ has 
a $3 \times 3$ principal submatrix with two positive eigenvalues. So, by 
Cauchy's interlacing principle, $L_r$ has at least two positive eigenvalues. 
Thus $L_r$ has exactly two positive and $n-2$ negative eigenvalues. In other 
words,
\begin{equation}
{\rm In} \,\,L_r = (2,0, n-2) \,\,\mbox{for}\,\, 3<r<4.   \label{eq21}
\end{equation}

At this stage note that the Theorem is completely proved for $n=4.$ Now let $n 
> 4,$ and consider the case $4 < r < 5.$ Arguing as before $\langle x, L_r x 
\rangle > 0$ for all $x \in \mathcal{H}_2.$ So $L_r$ has at least $n-2$ 
positive eigenvalues. It also has a $4 \times 4$ principal submatrix with two 
negative eigenvalues. Hence
\begin{equation}
{\rm In} \,\,L_r = (n-2, 0, 2) \,\,\,\mbox{for}\,\,\, 4 < r < 5.   
\label{eq22}
\end{equation}
The argument can be continued, introducing the space
\begin{eqnarray*}
 \mathcal{H}_3 &=& \left \{ x : \sum x_i = 0, \sum p_i x_i = 0, \sum p_i^2 x_i 
= 0   \right \} \\
&=&\left \{ x : Ex = 0, E Dx = 0, E D^2 x = 0    \right \}
\end{eqnarray*}
at the next stage. Using this we can prove statement (iii) for $5 < r < 6$ and 
$6 < r < 7.$ It is clear now how to complete the proof.

All parts of Theorem \ref{thm1.1} have now been established.$\hfill{\square}$

\vskip0.2in

We end this section with a few questions.

\begin{enumerate}
 \item Let $f(z)$ be the complex function defined as
$$f(z) = \det \left [\frac{p_i^z - p_j^z}{p_i - p_j} \right ]. $$
Our analysis has shown that $f$ has zeros at $z = 0, \pm 1, \pm 2, \ldots, \pm 
n-1;$ these zeros have multiplicities $n, n-1, \ldots, 1,$ respectively; and 
these are the only real zeros of $f.$ It might be of interest to find what 
other zeros $f$ has in the complex plane.

\item When $n=3,$ calculations show that
$$\det L_3 = - (p_1 - p_2)^2 (p_1 - p_3)^2 (p_2 - p_3)^2,$$
and
\begin{eqnarray*}
\det L_4 &=& - 2 (p_1 - p_2)^2 (p_1 - p_3)^2 (p_2 - p_3)^2 \\
&& \left \{(p_1 + p_2 + p_3) (p_1 p_2 + p_1 p_3 + p_2 p_3) + p_1 p_2 p_3 \right 
\}. 
\end{eqnarray*}
It might be of interest to find formulas for the determinants of the matrices 
$L_m$ for integers $m.$

\item Two of the authors have studied the matrix $P_r = \left [ (p_i + p_j)^r 
\right ]$ in \cite{bj}. It turns out that ${\rm In} \,\, P_r = {\rm In} 
\,\,L_{r+1}$ for all $r > 0.$ Why should this be so, and are there other 
interesting connections between these two matrix families?
\end{enumerate}


 \vskip0.1in
\noindent {\bf{Acknowledgements. }}The work of R. Bhatia is supported by a J. C. Bose National Fellowship, of S. Friedland by the NSF grant DMS-1216393, and of T. Jain by a SERB Women Excellence Award. The authors thank John Holbrook, Roger Horn, Olga Holtz and Lek-Heng Lim for illuminating discussions. The completion of this work was facilitated by the workshop ``Positivity, graphical models and modeling of complex multivariate dependencies'' at the American Institute of Mathematics in October 2014. The authors thank the organisers and the participants of this workshop.



\vskip0.2in
\begin{thebibliography}{99}
\itemsep10pt

\bibitem{rbh} R. Bhatia, {\it Matrix Analysis,} Springer, 1997.

\bibitem{rbh1} R. Bhatia, {\it Positive Definite Matrices,} Princeton
Univ. Press, 2007.

\bibitem{rb}
R. Bhatia, {\it First and second order perturbation bounds for the operator absolute value}, Linear Algebra Appl., 208 (1994) 367-376. 
\bibitem{bhol} R. Bhatia and J. A. Holbrook, {\it Fr\'echet derivatives of
the power function,} Indiana Univ. Math. J., 49(2000) 1155-1173.

\bibitem{bj} R. Bhatia and T. Jain, {\it Inertia of the matrix $\left [
(p_i+p_j)^r\right ],$} to appear in J. Spectral Theory.

\bibitem{bp}
R. Bhatia and K. R. Parthasarathy, {\it Positive definite functions and operator inequalities}, Bull. London Math. Soc., 32 (2000) 214-228.

\bibitem{rbsano} R. Bhatia and T. Sano, {\it Loewner matrices and operator
convexity,} Math. Ann., 344 (2009) 703-716.


\bibitem{bs}
R. Bhatia and K. B. Sinha, {\it Variation of real powers of positive operators,} Indiana Univ. Math. J., 43 (1994) 913-925.


\bibitem{hk1}
F. Hiai and H. Kosaki, {\it Comparison of various means for operators,} J. Funct. Anal., 163 (1999) 300-323.

\bibitem{hk2}
F. Hiai and H. Kosaki, {\it Means for matrices and comparison of their norms}, Indiana Univ. Math. J., 48 (1999) 899-936.

\bibitem{hk3}
F. Hiai and H. Kosaki, {\it Means of Hilbert Space Operators}, Springer 2003.

\bibitem{hn}
R. A. Horn, {\it The Hadamard product. Matrix theory and applications}, Proc. Sympos. Appl. Math., 40 (1989) 87-169.

\bibitem{hj}
R. A. Horn and C. R. Johnson, {\it Matrix Analysis}, Second edition, Cambridge Univ. Press, 2013.

\bibitem{hk4}
H. Kosaki, {\it Arithmetic-geometric mean and related inequalities for operators}, J. Funct. Anal., 156 (1998) 429-451.

\bibitem{k}
S. Karlin, {\it Total Positivity,} Stanford University Press, 1968.

 
\bibitem{gg} G. P\'olya and G. Szeg\"o, {\it Problems and Theorems in Analysis,}
Volume II, 4th ed., Springer, 1971.




\end{thebibliography}
\end{document}